\documentclass{svproc}
\usepackage{bm}
\usepackage{amsmath}
\usepackage{amsfonts}
\usepackage{mathtools}
\usepackage{enumitem} 
\usepackage{algorithmic} 
\usepackage{algorithm}
\usepackage{subfigure}

\usepackage{listings}
\newcommand{\beq}{\begin{equation}}
\newcommand{\eeq}{\end{equation}}

\def\Pi{\mathtt {Pi}}
\def\H{\mathtt {H}}

\def\h{\mathtt {h}}
\def\hedos{\mathtt {h^{2}}}
\def\Dhedos{\mathtt {Dh^{2}}}

\def\SNu{\mathtt {SN_1}}
\def\SNd{\mathtt {SN_2}}

\def\Puno{\mathtt {P^{1}}}
\def\Pdos{\mathtt {P^{2}}}
\def\H{\mathtt {H}}
\def\DH{\mathtt {DH}}

\def\Htres{\mathtt {H^{3}}}

\def\He{\mathtt {He}}
\def\DHe{\mathtt {DHe}}

\def\Hent{\mathtt {He^{34}}}

\def\T{\mathtt {T}}

\def\DH{\mathtt {DH}}

\def\DZ{\mathtt {DZ}}

\def\CU{\mathtt {CU}}
\def\TB{\mathtt {TB}}

\def\DHeuno{\mathtt {DHo^{E1}}}
\def\DHedos{\mathtt {DHo^{E2}}}

\def\DHeuno{\mathtt {DHe^{1}}}

\def\DHedos{\mathtt {DHe^{2}}}

\def\TPp{\mathtt {TP_p}}
\def\TPs{\mathtt {TP_s}}

\usepackage{url}

\begin{document}

\mainmatter              

\noindent
	This is a preprint of the following chapter: \newline
	Antonio Algaba, M. Cinta Dom\'inguez-Moreno, Manuel Merino,
	and Alejandro J. Rodr\'iguez-Luis,
	A degenerate Takens--Bogdanov bifurcation in a normal form of Lorenz's equations,
	published in 
	Advances in Nonlinear Dynamics
	Proceedings of the Second International Nonlinear Dynamics Conference
	(NODYCON 2021), Volume 1,
	edited by 
	W. Lacarbonara, B. Balachandran, M.J. Leamy, J. Ma, J.A. Tenreiro Machado, G. Stepan, 2022,
	Springer, Cham, 
	reproduced with permission of 
	Springer Nature Switzerland AG.
	The
	final authenticated version is available online at: http://dx.doi.org/10.1007/978-3-030-81162-4$\_$60. 

%
\title{A degenerate Takens--Bogdanov bifurcation in a normal form of Lorenz's equations}
\titlerunning{A degenerate double-zero bifurcation}  
%
\author{Antonio Algaba\inst{1} \and M. Cinta Dom\'inguez-Moreno \inst{1} \and Manuel Merino\inst{1} \and \\ 
Alejandro J. Rodr\'iguez-Luis\inst{2}}
\authorrunning{Antonio Algaba et al.} 
%
\tocauthor{Ivar Ekeland, Roger Temam, Jeffrey Dean, David Grove,
Craig Chambers, Kim B. Bruce, and Elisa Bertino}
\institute{
Departamento de Ciencias Integradas, Centro de Estudios Avanzados en F\'isica, Matem\'atica y Computaci\'on, 
Universidad de Huelva, 21071 Huelva, Spain\\
\email{algaba@uhu.es, mcinta.dominguez@dmat.uhu.es, merino@uhu.es}
\and
Departamento de Matem\'atica Aplicada II, E.T.S. Ingenieros, Universidad de Sevilla, 41092 Sevilla, Spain\\
\email{ajrluis@us.es}
}
\maketitle              

\begin{abstract}
In this work we consider an unfolding of a normal form of the Lorenz system near a triple-zero singularity. 
We are interested in the analysis of 
double-zero bifurcations 
emerging from that singularity.
Their local study provide partial results 
that are extended by means of numerical continuation methods. 
Specifically, a curve of heteroclinic connections is detected. 
It has a degenerate point from which infinitely many homoclinic connections emerge. 
In this way, we can partially understand the dynamics near the triple-zero singularity.
\keywords{double-zero, Takens--Bogdanov, Lorenz system, heteroclinic}
\end{abstract}

\section{Introduction}
In the world of dynamical systems
the first known and most famous chaotic one is the Lorenz system \cite{Lo:63}
\begin{equation}
\label{Lorenz}
\dot x = \sigma (y - x), \ \ \ \
\dot y = \rho x - y - xz, \ \ \ \
\dot z = - b z + xy, \ \ \ \ \ \ \ \ 
\sigma, \rho, b \in \mathbb{R}. 
\end{equation}
The study of its riveting and intricate dynamical behavior has been carried out in multitude of works
(see the recent papers 
\cite{DoKrOs:11,BaBlSe:11,BaShSh:12,AlFeMeRo:14b,AlFeMeRo:15,CrKrOs:15a,DoKrOs:15,%
AlDoMeRo:15,AlDoMeRo:16,AlGaMeRo:16,AlMeRo:16,CrKrOs:17,Os:18,AlDoMeRo:18}
and references therein).

One way to obtain important information on the organizing centers of the dynamics in system (\ref{Lorenz})
is by means of the study of local bifurcations of equilibria. 
Whereas the Hopf and Takens-Bogdanov bifurcations have been fully studied in Lorenz system
\cite{AlDoMeRo:15,AlDoMeRo:16}, 
the Hopf-pitchfork bifurcation (a pair of imaginary eigenvalues and the third one zero; it occurs
when $\sigma = -1$,  $b = 0$, $\rho > 1$)
and the triple-zero bifurcation (a triple-zero eigenvalue; it arises if $\sigma = -1$, $b = 0$, $\rho = 1$)
cannot be analyzed by the standard procedures because non-isolated equilibria appear when $b=0$.

To avoid this problem, in this paper we consider a three-parameter unfolding, that is close to the normal form
of the triple-zero bifurcation exhibited by Lorenz system, given by \cite{AlDoMeRo:20}
\begin{equation}
\label{unfoldingFNLorenz1}
\dot x = y, \ \ \ \
\dot y = \varepsilon_1 x + \varepsilon_2 y + A xz + B yz, \ \ \ \
\dot z = \varepsilon_3 z + C x^2 + D z^2,
\end{equation}
where $\varepsilon_1, \varepsilon_2, \varepsilon_3 \approx 0$ 
and $A,B,C,D$ are real parameters.
System (\ref{unfoldingFNLorenz1}) exhibits a triple-zero bifurcation when 
$\varepsilon_1=\varepsilon_2=\varepsilon_3 = 0$. These equations are also invariant under the change
$(x,y,z) \rightarrow (-x,-y,z)$.
We remark that several systems studied in the literature appear 
as particular cases of (\ref{unfoldingFNLorenz1}) for certain parameter choices
\cite{Shim:80,Shil:93,Ru:93,Liu:04,Me:08,KoRo:04}.
Without loss of generality \cite{AlDoMeRo:20}, we can take $A=-1, C=1$
\begin{equation}
\label{unfoldingFNLorenz}
\dot x  =  y, \ \ \ \ 
\dot y  =  \varepsilon_{1}x+\varepsilon_{2}y-xz+Byz, \ \ \ \
\dot z  = \varepsilon_{3}z+x^2+Dz^2.
\end{equation}

System (\ref{unfoldingFNLorenz}) can have up to four equilibria, namely
$E_1=(0,0,0)$,
$E_2=\left(0,0,{-\varepsilon_3}/{D}\right)$ if  $D\neq 0$ 
and
$E_{3,4} = (\pm\sqrt{-\varepsilon_1(\varepsilon_3+ D\varepsilon_1)},0,\varepsilon_1 )$
if  $\varepsilon_1(\varepsilon_3+ D\varepsilon_1) < 0$. 
Note that $E_1$ and $E_2$ are placed on the $z$-axis, which is an invariant set.

The characteristic polynomial of the Jacobian matrix of system (\ref{unfoldingFNLorenz}) at the origin $E_1$ 
is given by
$P(\lambda) = \lambda^3+p_1\lambda^2+p_2\lambda+p_3$, 
where
$
p_1=-(\varepsilon_2+\varepsilon_3), \
p_2=\varepsilon_2\varepsilon_3-\varepsilon_1, \
p_3=\varepsilon_1\varepsilon_3.
$
Therefore, the origin exhibits the following bifurcations:
(a)
A pitchfork bifurcation when 
$\varepsilon_1=0$, 
$\varepsilon_2 \neq 0$, 
$\varepsilon_3 \neq 0$.
The nontrivial equilibria $E_3$ and $E_4$ appear when 
$-\varepsilon_1(\varepsilon_3+D\varepsilon_1)>0$.
(b)
A transcritical bifurcation (involving also $E_2$) when 
$\varepsilon_3=0$, 
$\varepsilon_1 \neq 0$, 
$\varepsilon_2 \neq 0$,
$D \neq 0$.
(c)
A Hopf bifurcation when
$\varepsilon_1 < 0$, 
$\varepsilon_2 = 0$, 
$\varepsilon_3 \neq 0$.
(d)
A Takens--Bogdanov bifurcation 
(nondiagonalizable double-zero eigenvalue)
when
$\varepsilon_1=0$, 
$\varepsilon_2=0$, 
$\varepsilon_3 \neq 0$.
It is of homoclinic type when $\varepsilon_3<0$ and of heteroclinic type if $\varepsilon_3>0$.
(e) 
A Hopf-zero bifurcation when
$\varepsilon_2=0$, 
$\varepsilon_3=0$, 
$\varepsilon_1<0$.
(f)
A double-zero bifurcation (a diagonalizable double-zero eigenvalue)
when
$\varepsilon_1=0$, 
$\varepsilon_3=0$, 
$\varepsilon_2 \neq 0$ (see \cite{AlDoMeRo:20}).
(g)
A triple-zero bifurcation when $\varepsilon_1=\varepsilon_2=\varepsilon_3=0$.

We show now that all the information on the equilibrium $E_2$ can be obtained from the
analysis of $E_1$.
To study $E_2$
we translate it to the origin by means of the change
$
x={\tilde x}$,  
$y={\tilde y}$, 
$z={\tilde z}-\varepsilon_3/D,
$
that transforms system (\ref{unfoldingFNLorenz}) into
\begin{equation}
\label{HopfE2}
\dot {\tilde x}  =  {\tilde y},  \ \ \
\dot {\tilde y}  =  (\varepsilon_1 + \frac{1}{D} \varepsilon_3) {\tilde x}
                         + (\varepsilon_2-\frac{B}{D} \varepsilon_3) {\tilde y} - {\tilde x} {\tilde z} + B {\tilde y} {\tilde z}, \ \ \
\dot {\tilde z}  =  -\varepsilon_3 {\tilde z} + {\tilde x}^2 + D{\tilde z}^2,
\end{equation}
with $\varepsilon_3, D \neq 0$.

Since system (\ref{unfoldingFNLorenz}) is symmetric to the change
\begin{equation}
\label{Cambiodeva}
\hspace{-0.03truecm}
(x,y,z,t,\varepsilon_1, \varepsilon_2, \varepsilon_3,B,D)  
\hspace{-0.05truecm}
\rightarrow
\hspace{-0.05truecm}
\left(    
    \hspace{-0.05truecm}
    x, y, z-\displaystyle \frac{\varepsilon_3}{D},t,
     \varepsilon_1+ \displaystyle\frac{\varepsilon_3}{D}, 
     \varepsilon_2 -\displaystyle\frac{B\,\varepsilon_3}{D}, -\varepsilon_3,B,D 
     \hspace{-0.05truecm}
\right) 
     \hspace{-0.05truecm} ,     
\end{equation}
it is direct
to obtain the stability and bifurcations of $E_2$ from the stability and bifurcations of $E_1$.
Thus, it is enough to study the bifurcations exhibited by $E_1$.

\section{Double-zero bifurcations of the origin}
\label{sec:TB}

The Takens--Bogdanov bifurcation of $E_1$ occurs when $p_{2}=p_{3}=0,$ $p_{1}\neq0$,
that is, when $\varepsilon_{1}=\varepsilon_{2}=0$, $\varepsilon_{3}\neq0$.
For these values the Jacobian matrix has a nondiagonalizable double-zero and a nonzero eigenvalue. 
Thus, system (\ref{unfoldingFNLorenz}) reads
\begin{equation}
\label{sys01}
\dot x  =  y,  \ \ \ \
\dot y  =  -xz+Byz,   \ \ \ \
\dot z  =  \varepsilon_{3}z+x^2+Dz^2.
\end{equation}

By means of the second-order approximation to the center manifold
we obtain the third-order reduced system and compute its normal form
\begin{equation}
\label{sys02}
\dot x  =  y, \ \ \ \  
\dot y  =  a_{3}x^3+b_{3}x^2y, \ \ \ \ 
{\rm with} \ \ \ \ a_{3} = 1/\varepsilon_3, 
b_{3} = (2-\varepsilon_{3}B)/\varepsilon_{3}^2 .
\end{equation}

Its unfolding is given by
\begin{equation}
\label{sys03}
\dot x  =  y, \ \ \ \
\dot y  =  \varepsilon_1 x+\varepsilon_2 y+a_{3}x^3+b_{3}x^2y.
\end{equation}
By means of the rescaling
\begin{equation}
\label{scal01}
x \to \frac{|\varepsilon_{3}|\sqrt{|\varepsilon_{3}|}}{2-B\varepsilon_{3}}  \, \bar{x},
\quad 
y \to \frac{\varepsilon_{3}^2\sqrt{|\varepsilon_{3}|}}{(2+B\varepsilon_{3}^2)}  \, \bar{y},
\quad 
t \to \frac{2-B\varepsilon_{3}}{|\varepsilon_{3}|}  \, \tau,
\end{equation}
system (\ref{sys03}) is transformed into
\begin{equation}
\label{sys04}
\dot {\bar{x}}  =  \bar{y}, \ \ \ \
\dot {\bar{y}}  =  \frac{\varepsilon_{1}(2-B\varepsilon_{3})^2}{\varepsilon_{3}^2}\bar{x}
+\frac{\varepsilon_{2}(2-B\varepsilon_{3})}{|\varepsilon_{3}|}\bar{y}
+sgn(\varepsilon_{3})\bar{x}^3-\bar{x}^2\bar{y}.
\end{equation}

Thus, the Takens--Bogdanov bifurcation is of heteroclinic case when $\varepsilon_{3}>0$
and of homoclinic case if $\varepsilon_{3}<0$.
Note that if $\varepsilon_{3}=0$ a triple-zero bifurcation is present. 
Moreover, when
$(\varepsilon_{1},\varepsilon_{2},\varepsilon_{3})=(0,0,2/B)$, 
with $B\neq0$, a degenerate Takens--Bogdanov bifurcation occurs. This case will be analyzed below.

In the homoclinic case, when $\varepsilon_{3}<0$,  the equilibria of system (\ref{sys04}) are
$(0,0)$ and 
$(\pm (2-B\varepsilon_{3})\sqrt{\varepsilon_1}/\varepsilon_3,0)$, 
$\varepsilon_{3} \neq 2/B$, $\varepsilon_{1}>0$.
From the Takens-Bogdanov singularity the curves corresponding to the following bifurcations emerge 
\cite{GuHo:83,Wi:03,Ku:04}:
(a)
A pitchfork bifurcation of the origin is present when  $\varepsilon_{1}=0$. 
(b)
A subcritical Hopf bifurcation of the origin if  
$\varepsilon_{2}=0$, $\varepsilon_{1}<0$. 
(c)
A supercritical Hopf bifurcation of the nontrivial equilibria when
$\varepsilon_{2} \approx (2-B\varepsilon_{3}) \varepsilon_1/\varepsilon_3$, 
for
$\varepsilon_{1}>0$. 
(d)
A homoclinic connection to the origin for 
$\varepsilon_{2} \approx 4(2-B\varepsilon_{3}) \varepsilon_{1} /(5\varepsilon_3)$, 
with $\varepsilon_{1}>0$. 
Since the third eigenvalue ($\varepsilon_{3}<0$) determines the behavior outside the center manifold, 
these homoclinic connections are attractive.
(e)
A saddle-node bifurcation of symmetric periodic orbits when
$\varepsilon_{2} \approx c \, \varepsilon_{1}(2-B\varepsilon_{3})/\varepsilon_{3}$, 
where $\varepsilon_{1}>0$, $c\approx 0.752$.

In the heteroclinic case, when $\varepsilon_{3}>0$,  the equilibria of system (\ref{sys04}) are
$(0,0)$ and 
$(\pm (2-B\varepsilon_{3})\sqrt{-\varepsilon_{1}} / \varepsilon_{3},0)$,
$\varepsilon_{3}\neq 2/B$, $\varepsilon_{1}<0$.
The following curves are present \cite{GuHo:83,Wi:03,Ku:04}:
(a)
A pitchfork bifurcation of the origin for 
(b)
A subcritical Hopf bifurcation of the origin for 
(c) 
A heteroclinic connection to nontrivial equilibria if
$\varepsilon_{2} \approx (2-B\varepsilon_3)\varepsilon_1/(5\varepsilon_3)$,  $\varepsilon_{1}<0$. 
As the third eigenvalue ($\varepsilon_{3}>0$) determines the behavior outside the center manifold, 
these heteroclinic connections are repulsive.

\subsection{Codimension-three degeneracy}

The normal form coefficient $b_{3}$, given in  (\ref{sys02}), vanishes when
$\varepsilon_{3}={2/B}$, $B\neq0$. 

Thus, considering the fourth-order approximation to the center manifold
we obtain the fifth-order reduced system on the center manifold 
\begin{equation}
\label{sys07}
\dot x  =  y, \ \ \  
\dot y  = a_{3}x^3+a_{5}x^5+b_{5}x^4y. 
\end{equation}
Multiplying system (\ref{sys07}) by $(1-(a_5/a_3)x^2)$, 
the $x^5$-term can be eliminated
\begin{equation}
\label{sys08}
\dot x  =  y+\frac{-a_{5}}{a_{3}}x^2y, \ \ \ 
\dot y  =  a_{3}x^3+b_{5}x^4y,
\ \  {\rm with} \ \ 
a_{3}=\frac{1}{\varepsilon_{3}}, \ b_{5}=\frac{-B^4}{8}(5B+3D).
\end{equation}

An unfolding is given by 
\begin{equation}
\label{sys09}
\dot x  =  y,\ \ \ \
\dot y  =  \mu_{1}x+\mu_{2}y+a_{3}x^3+\mu_{3}x^2y+b_{5}x^4y,
\end{equation}
with
$\mu_{1}=\varepsilon_{1}$, 
$\mu_{2}=\varepsilon_{2}$, 
$\mu_{3}=(2-B\varepsilon_{3})/\varepsilon_{3}^2$.

To analyze this degeneracy we use the rescaling
$$
x \to \sqrt[6]{\frac{|a_{3}|}{b_{5}^2}} \, u, \ \ \ 
y \to b_{5}\sqrt[6]{\frac{|a_{3}|^5}{b_{5}^{10}}} \, v, \ \ \
t \to \sqrt[3]{ \frac{b_5}{|a_3|^2} } \, \tau,
$$
and then system (\ref{sys09}) is transformed into
\begin{equation}
\label{sys10}
\dot u  =  v,\ \ \ \
\dot y  =  \tilde{\mu}_{1}u+\tilde{\mu}_{2}v+sgn(a_{3})u^3+\tilde{\mu}_{3}u^2v+u^4v,
\end{equation}
where
$$
\tilde{\mu}_{1}
\hspace{-0.1truecm} = \hspace{-0.1truecm}
\frac{\varepsilon_{1}}{4} \hspace{-0.05truecm} \sqrt[3]{\varepsilon_{3}^4B^8(5B+3D)^2},
\
\tilde{\mu}_{2}
\hspace{-0.1truecm} = \hspace{-0.1truecm}
\frac{-\varepsilon_{2}}{2}\sqrt[3]{-B^4\varepsilon_{3}^2(5B+3D)},
\
\tilde{\mu}_{3}
\hspace{-0.1truecm} = \hspace{-0.1truecm}
\frac{-2(2-B\varepsilon_{3})\sqrt[3]{\varepsilon_{3}}}{\sqrt[3]{B^4(5B+3D)}} .
$$

The analysis of this three-parameter family can be found in \cite{LiRo:90,RoFrPo:91}.
In the parameter space several codimension-two bifurcation curves emerge from 
the point $(0,0,2/B)$.

In the homoclinic case ($sgn(a_{3})<0$, i.e., $\varepsilon_{3}<0$):
\begin{itemize}
\item 
A nondegenerate Takens-Bogdanov for 
$\tilde{\mu}_{1}=\tilde{\mu}_{2}=0, \tilde{\mu}_{3} \neq 0$, 
i.e., 
$\varepsilon_{1}=\varepsilon_{2}=0$, 
$\varepsilon_{3} \neq 2/B$.

\item 
A degenerate Hopf bifurcation of the origin for
$\tilde{\mu}_{2}=\tilde{\mu}_{3}=0, \tilde{\mu}_{1}<0$, 
i.e.,
$\varepsilon_{2}=0$, 
$\varepsilon_{3}=\frac{2}{B}$, 
$\varepsilon_{1}<0$. 

\item
A degenerate Hopf bifurcation of the nontrivial equilibria when
$\tilde{\mu}_{2}=-\tilde{\mu}_{1}^2$, $\tilde{\mu}_{3}=0$, i.e.,
$\varepsilon_{3}=\frac{2}{B}$,
$\varepsilon_{2}=\frac{-1}{2}\varepsilon_{1}^2B^2(5B+3D)$, 
$\varepsilon_{1}>0$.

\item 
A degenerate (zero-trace) homoclinic connection to the origin when 
$\tilde{\mu}_{2}=0$, 
$\tilde{\mu}_{3}=\frac{-8\tilde{\mu_{1}}}{7}$, 
$\tilde{\mu}_{1}>0$,
i.e., 
$\varepsilon_{2}=0$, 
$\varepsilon_{1}=\frac{7(2-B\varepsilon_{3})}{B^4\varepsilon_{3}(5B+3D)}>0$.

\item 
A cusp of saddle-node bifurcations of periodic orbits when 
$\tilde{\mu}_{2}=c_{3}\tilde{\mu}_{1}^2$, 
$\tilde{\mu}_{3}=c_{4}\tilde{\mu}_{1}$, 
$\tilde{\mu}_{1}>0$, 
with $c_{3}\approx 1.5713$, $c_{4}\approx-3.3484$.
That is, for  
$\varepsilon_{2}=\frac{2c_{3}(2-B\varepsilon_{3})^2}{c_{4}^2B^2(5B+3D)}$, 
$\varepsilon_{1}=\frac{-8(2-B\varepsilon_{3})}{c_{4}\varepsilon_{3}B^4(5B+3D)}>0$.

\end{itemize}

In the heteroclinic case ($sgn(a_{3})>0$, i.e., $\varepsilon_{3}>0$):
\begin{itemize}
\item 
A nondegenerate Takens-Bogadnov bifurcation when 
$\tilde{\mu}_{1}=\tilde{\mu}_{2}=0$, 
$\tilde{\mu}_{3} \neq0$, 
i.e., 
$\varepsilon_{1}=\varepsilon_{2}=0$, 
$\varepsilon_{3} \neq 2/B$.

\item 
A degenerate Hopf bifurcation of the origin for
$\tilde{\mu}_{2}=\tilde{\mu}_{3}=0$, 
$\tilde{\mu}_{1}<0$, 
i.e., 
$\varepsilon_{2}=0$, 
$\varepsilon_{3}=2/B$, 
$\varepsilon_{1}<0$.

\item 
A degenerate (zero-trace) heteroclinic connection when 
$\tilde{\mu}_{2}=0$, 
$\tilde{\mu}_{3} \approx \frac{3}{7}\tilde{\mu}_{1}$, 
$\tilde{\mu}_{1}<0$, 
i.e.,  
$\varepsilon_{2}=0$, 
$\varepsilon_{1} \approx -\frac{7}{3}(2-B\varepsilon_{3})^3$.

\end{itemize}

We end this section with two useful comments for Sect. \ref{Theonumer}.
On the one hand, the information on the Takens--Bogdanov bifurcation 
exhibited by $E_2$ can be easily obtained from the above results by using the symmetry
(\ref{Cambiodeva}): it is of homoclinic-type if $\varepsilon_3>0$ and of heteroclinic-type when $\varepsilon_3<0$.
On the other hand, the (diagonalizable) double-zero bifurcation exhibited by the origin has been analyzed in 
\cite{AlDoMeRo:20}. There, the following expression is obtained for the curve of heteroclinic connections between
$E_1$ and $E_2$ (in the nondegenerate case $B \neq 2 \Delta - \Delta^3 \varepsilon_1$)
$$
\varepsilon_1  
\approx 
-\frac{1}{2D} \, \varepsilon_3 + \frac{a^2(-2\Delta+B)}{4D(3a+2)(-\Delta^2\varepsilon_1+1)} \,  \varepsilon_3^2,
\ \ {\rm where} \ \ \Delta=\frac{1}{\varepsilon_2}, \ 
a = \frac{\Delta^3 \varepsilon_1 - \Delta}{D}>0.
$$

\section{Numerical study}
\label{Theonumer}

With the information provided in Sect. \ref{sec:TB},
we are going to perform a numerical study 
of system (\ref{unfoldingFNLorenz}) with the continuation software AUTO \cite{DoOl:12}.
Specifically, we will draw bifurcation sets in the $(\varepsilon_1,\varepsilon_3)$-parameter plane, 
in a neighbourhood of the degeneracies $\DZ$ (diagonalizable double-zero) and $\TB$ (Takens--Bogdanov).
We will fix $\varepsilon_2=-1$, $B=-0.1<0$ and $D=0.01>0$ according to the values used in \cite{KoRo:04}.

\begin{figure}[t]
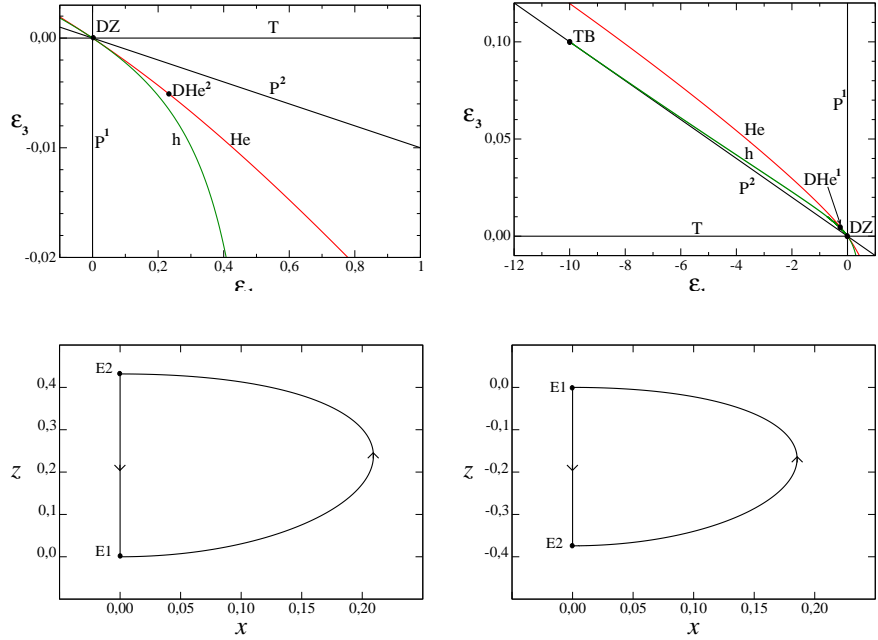

\begin{center}
{\bf (a)} \hspace{5.5truecm} {\bf (b)}
\end{center}
\vspace{-0.6cm}
\begin{center}
\includegraphics[angle=0,width=5.5cm]{Figura1a.eps}
\hspace{0.3cm}
\includegraphics[angle=0,width=5.5cm]{Figura1b.eps}
\end{center}
\vspace{-0.7cm}
\begin{center}
{\bf (c)} \hspace{5.5truecm} {\bf (d)}
\end{center}
\vspace{-0.6cm}
\begin{center}
\includegraphics[angle=0,width=5.5cm]{Figura1c.eps}
\hspace{0.3cm}
\includegraphics[angle=0,width=5.5cm]{Figura1d.eps}
\end{center}
\vspace{-0.5cm}
\caption{
\label{Figure1} 
For $\varepsilon_2=-1, B=-0.1, D=0.01$ partial bifurcation set in a neighborhood 
of the point $\DZ$: 
(a) in the fourth quadrant;
(b) in the second quadrant. 
(c) Projection onto the $(x,z)$-plane of the heteroclinic cycle $\He$ of panel (a) that exists 
     when $(\varepsilon_1,\varepsilon_3) \approx (0.2, -0.0043175)$. 
(d) Projection onto the $(x,z)$-plane of the heteroclinic cycle $\He$ of panel (b) that exists 
     when $(\varepsilon_1,\varepsilon_3) \approx (-0.2, 0.0037414)$.
}
\end{figure}

First we draw, in Figs. \ref{Figure1}(a) and \ref{Figure1}(b), 
partial bifurcation sets with the bifurcation curves related to the degeneracies
$\DZ$ (of the origin) and $\TB$ (of $E_2$) in the fourth and second quadrants, respectively.
We can observe the curves $\h$ (Hopf bifurcation of the equilibria $E_{3,4}$) and 
$\He$ (heteroclinic cycle between the equilibria $E_1$ and $E_2$).
Both curves emerge from the point where the double-zero degeneracy of the origin occurs, $\DZ = (0,0)$ 
in the $(\varepsilon_1,\varepsilon_3)$-parameter plane. 
Moreover, three straight lines intersect at the double-zero point $\DZ$, namely 
$\Puno$ (pitchfork bifurcation of the origin, $\varepsilon_1=0$), 
$\Pdos$ (pitchfork bifurcation of $E_2$, $\varepsilon_3=-0.01 \varepsilon_1$) and 
$\T$ (transcritical bifurcation between $E_1$ and $E_2$, $\varepsilon_3=0$).
Note that in these pitchfork bifurcations the equilibria $E_{3,4}$ emerge. 
Fixed a value of $\varepsilon_3$ in a neighborhood of the point $\DZ$, 
the periodic orbit emerged from $\h$ is atractive and, increasing the value of $\varepsilon_1$,  
it disappears in a heteroclinic cycle $\He$.

This cycle $\He$ is formed by two heteroclinic connections, 
one is structurally stable (since it goes from $E_1$ to $E_2$ on the invariant $z$-axis) 
and the other one is more relevant (because it is placed outside the $z$-axis). 
To show the differences that exist between the heteroclinic cycles of 
Figs. \ref{Figure1}(a) and \ref{Figure1}(b) we draw their projections onto the $(x,z)$-plane 
in Figs. \ref{Figure1}(c) and \ref{Figure1}(d), respectively.
Remark that throughout this work, with the aim of simplifying the notation, 
we will label the heteroclinic cycle (in fact, due to the symmetry, a pair of heteroclinic cycles exist) 
and the heteroclinic bifurcation with the same symbol, although they are two different objects. 
Also, when necessary, we will use superscripts to indicate the equilibria that are involved in a certain bifurcation 
or in its degeneration.

Such heteroclinic cycles are attractive when they emerge from $\DZ$
since the saddle quantities 
$
\delta_{E_1}=\left|max(\lambda_1,\lambda_3)/\lambda_2\right|
$
and
$
\delta_{E_2}=\left| \lambda^{*}_2/\lambda^{*}_3\right| 
$
satisfy 
$
\delta_{E_1} \delta_{E_2}>1
$.
Here we denote the eigenvalues of the Jacobian matrix at the origin $E_1$
as $\lambda_1,\lambda_3<0<\lambda_2$, where
$
\lambda_3 = \varepsilon_3
$,
$
\lambda_{2,1} = [\varepsilon_2 \pm \sqrt{\varepsilon_2^ 2 + 4 \varepsilon_1} \, ]/2 ,
$
and 
the eigenvalues of the Jacobian matrix at $E_2=(0,0,-100 \,\varepsilon_3)$
as $\lambda^{*}_1 < \lambda^{*}_2<0<\lambda^{*}_3$, with
$
\lambda^{*}_3= -\varepsilon_3
$, 
$\lambda^{*}_{2,1} = [ \varepsilon_2-\frac{B}{D}\varepsilon_3 \pm
\sqrt{\left(\varepsilon_2 -\frac{B}{D}\varepsilon_3\right)^ 2 +
4\left(\varepsilon_1+\frac{\varepsilon_3}{D}\right)} \, ] /2
$.

\begin{figure}[tb]
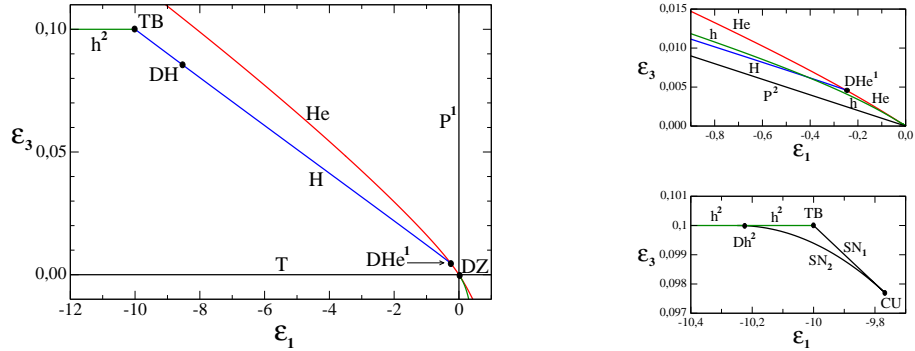

\begin{center}
\hspace{1.8truecm} {\bf (a)} \hspace{6truecm} {\bf (b)}
\end{center}
\vspace{-0.5cm}
\begin{center}
\includegraphics[angle=0,height=4.6cm]{Figura2a.eps}
\hspace{1.7cm}
\includegraphics[angle=0,height=4.6cm]{Figura2b.eps}
\end{center}
\vspace{-0.6cm}
\caption{
\label{Figure2} 
For $\varepsilon_2=-1, B=-0.1, D=0.01$ partial bifurcation set: 
(a) in the second quadrant; 
(b) zoom of panel (a) in the vicinity of the points $\DHeuno$ (upper panel) and $\TB$ (lower panel). 
}
\end{figure}

\begin{figure}[tb]
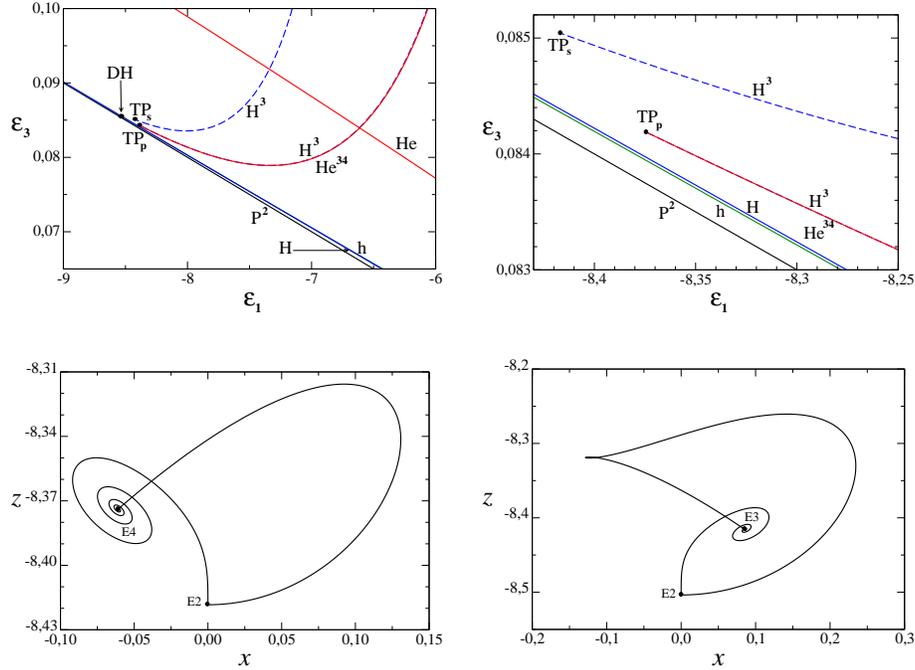

\begin{center}
{\bf (a)} \hspace{5.5truecm} {\bf (b)}
\end{center}
\vspace{-0.6cm}
\begin{center}
\includegraphics[angle=0,width=5.75cm]{Figura3a.eps}
\hspace{0.3cm}
\includegraphics[angle=0,width=5.75cm]{Figura3b.eps}
\end{center}
\vspace{-0.7cm}
\begin{center}
{\bf (c)} \hspace{5.5truecm} {\bf (d)}
\end{center}
\vspace{-0.6cm}
\begin{center}
\includegraphics[angle=0,width=5.75cm]{Figura3c.eps}
\hspace{0.3cm}
\includegraphics[angle=0,width=5.75cm]{Figura3d.eps}
\end{center}
\vspace{-0.5truecm}
\caption{
\label{Figure3} 
For $\varepsilon_2=-1$, $B=-0.1$, $D=0.01$:
(a) partial bifurcation set in the second quadrant.
(b) Zoom of panel (a). 
Projection onto the $(x,z)$-plane of T-point heteroclinic loops 
connecting $E_2$ and the equilibria $E_{3,4}$
(note that because of the symmetry a pair of the corresponding orbits exists): 
(c) principal T-point for $(\varepsilon_1,\varepsilon_3)\approx(-8.3738877, 0.0841835)$; 
(d) secondary T-point when $(\varepsilon_1,\varepsilon_3)\approx(-8.4159326, 0.0850368)$.
}
\end{figure}

\begin{figure}[tb]
\begin{center}
{\bf (a)} \hspace{5.5truecm} {\bf (b)}
\end{center}
\vspace{-0.77truecm}
\begin{center}
\includegraphics[angle=0,width=5.55cm]{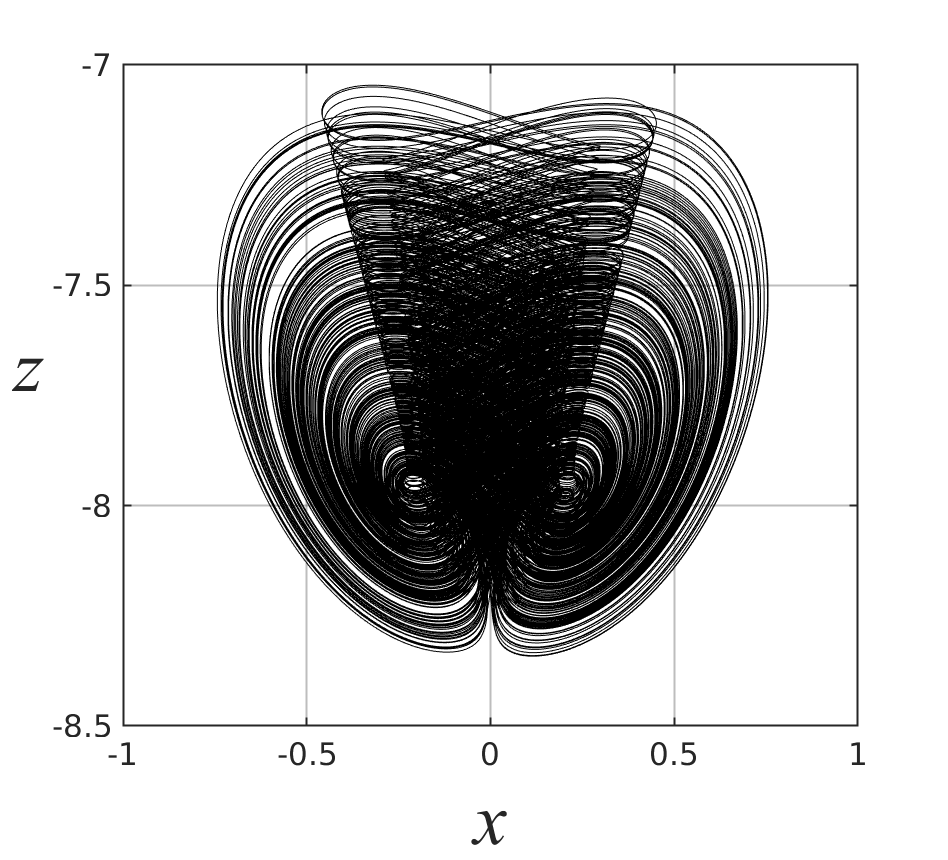}
\hspace{0.3cm}
\includegraphics[angle=0,width=5.55cm]{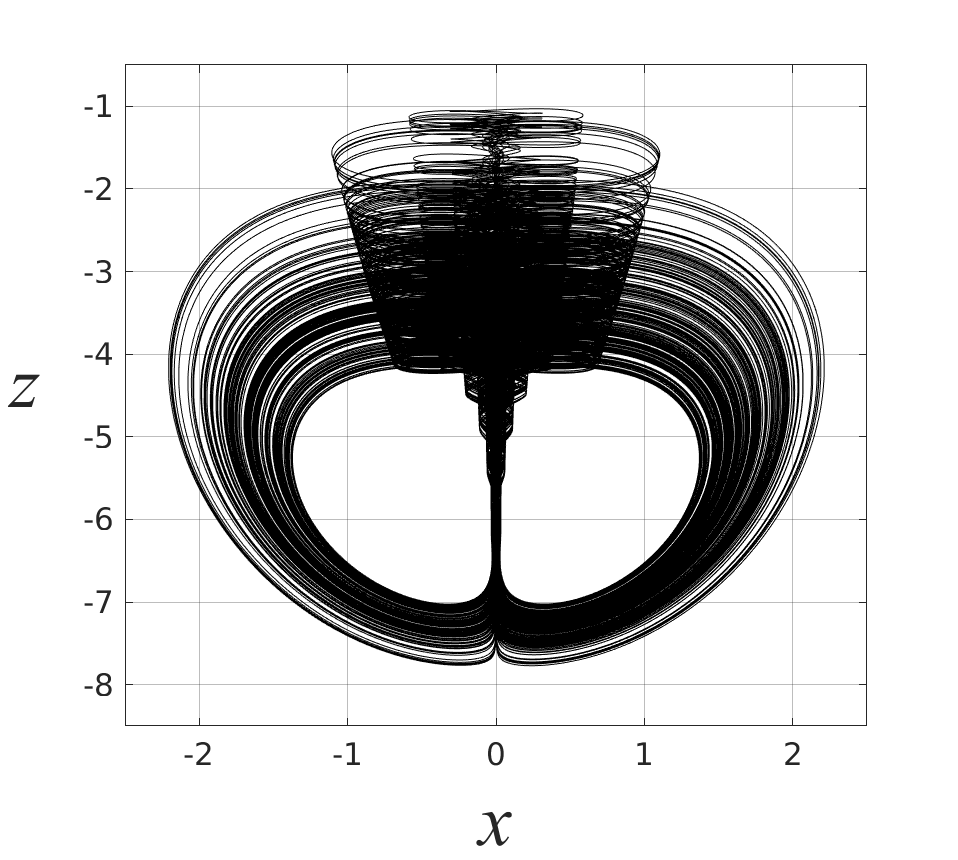}
\end{center}
\vspace{-0.7cm}
\caption{
\label{Figure4} 
For $\varepsilon_2=-1, \varepsilon_3=0.085, B=-0.1, D=0.01$ chaotic attractor when: 
(a) $\varepsilon_1=-8$, with initial conditions $(x_0,y_0,z_0)=(0,0.1,-8)$.
(b) $\varepsilon_1=-6.3$, with initial conditions $(x_0,y_0,z_0)=(0,0.1,-7)$.
}
\end{figure}

As seen in Figs. \ref{Figure1}(a) and \ref{Figure1}(b), 
on each of the curves $\He$, there is a point $\DHe$ where the heteroclinic cycle is degenerate.
In fact, when $\varepsilon_1>0$, the equilibrium $E_1$ is always a real saddle along the curve $\He$
but $E_2$, that is also a real saddle when it arises from $\DZ$, becomes a saddle-focus from the point
$\DHedos\approx (0.2328879,-0.0050898)$.
In the case of the branch located in the fourth quadrant, when
$\varepsilon_1<0$, the equilibrium $E_2$ is always a real saddle along the curve $\He$, 
whereas the equilibrium $E_1$, that emerges from $\DZ$ as a real saddle,
becomes a saddle-focus from $\DHeuno\approx (-0.25,0.004611)$.

The change in the configuration of $E_1$ (resp. $E_2$) on the curve of the heteroclinic connections $\He$
in the second (resp. fourth) quadrant implies the appearance of an infinity of bifurcation curves, 
which arise from the aforementioned point $\DHeuno$ (reps. $\DHedos$) on the $\He$ curve.

Specifically, in Figs. \ref{Figure2}(a) and \ref{Figure2}(b), we can see the curve 
$\H$ (of homoclinic connections to the equilibrium $E_2$) that emerges from the point $\DHeuno$.
In a vicinity of $\DHeuno$, a saddle periodic orbit emerges from $\H$ 
(since $\delta_{E2}<1$ at the points of such a neighborhood).
A degenerate point $\DH$ appears on $\H$ when 
$(\varepsilon_1,\varepsilon_3) \approx (-8.5339606, 0.0855365)$ 
because $\delta_{E2}=1$.
The curve $\H$ ends at the point $\TB=(-10,0.1)$ where $E_2$ undergoes a Takens-Bogdanov bifurcation.
This is in agreement with the theoretical results of Sect.~\ref{sec:TB} that guarantee, for these parameter values, 
that the Takens-Bogdanov bifurcation of $E_2$ is of homoclinic type.

In the second quadrant, the Hopf bifurcation $\h$ of the nontrivial equilibria $E_{3,4}$  is always supercritical 
and this curve connects the points $\DZ$ and $\TB$ (see Fig. \ref{Figure1}(b)).
In addition to the homoclinic connection curve $\H$, other curves emerge from $\TB$,
namely a curve $\hedos$ of Hopf bifurcation of $E_2$ ($\varepsilon_3=0.1$ when $\varepsilon_1<-10$)
and a curve $\SNu$ of saddle-node bifurcation of periodic orbits. 
The Hopf bifurcation $\hedos$ is subcritical when it emerges from $\TB$ and it becomes supercritical
because a degeneracy occurs at $\Dhedos\approx(-10.2487498,0.1)$ 
as the first Lyapunov coefficient vanishes.
A new curve $\SNd$ of saddle-node bifurcations of symmetric periodic orbits emerges from $\Dhedos$. 
A cusp bifurcation of periodic orbits $\CU\approx(-9.7715876,0.0977398)$ occurs when $\SNu$ and $\SNd$
collapse (see lower panel of Fig. \ref{Figure2}(b)).
The cusp $\CU$ is the first one of an infinite sequence of cusps that accumulate to the point $\DH$ \cite{Wi:03}.

Close to the degeneracy $\DH$,  other bifurcation curves
can be seen in Figs. \ref{Figure3}(a) and \ref{Figure3}(b). 
Specifically from the point
$\TPp \approx (-8.3738877, 0.0841835)$ (where a T-point heteroclinic loop between $E_2$ and $E_{3,4}$ 
exists and whose projection in the $(x,z)$-plane can be seen in Fig. \ref{Figure3}(c)) three curves of global
connections arise, namely
$\Hent$ (heteroclinic connection between $E_3$ and $E_4$), 
$\Htres$ (homoclinic connection to $E_{3,4}$) 
and a curve of homoclinic connections to $E_2$ (not included in Fig. \ref {Figure3}) 
that ends at $\DHeuno$
(see Figs. \ref{Figure2}(a)  and \ref{Figure2}(b))
\cite{GlSp:86,FeFrRo:02,FeFrRo:08,AlFeMeRo:15}.
The curve $\Htres$ ends in a secondary T-point 
between $E_2$ and $E_ {3,4}$ for
$\TPs \approx (-8.4159326, 0.0850368)$ whose projection on the $(x,z)$-plane appears in
Fig. \ref{Figure3}(d).

The presence in the second quadrant of the degenerations analyzed 
(and the bifurcation curves that arise from them)
determines the existence of regions where chaotic attractors exist 
(see Figs. \ref{Figure4}(a) and \ref{Figure4}(b)). 
We note that these attractors, obtained for $\varepsilon_3 = 0.085$, are structurally stable when 
the value of $\varepsilon_1$ is varied. Moreover, as observed in these figures, 
their size increases when $\varepsilon_1$  augments along the interval $[-8, -6.3]$.

\section{Conclusions}
\label{sec:conclusions}

In this work we consider an unfolding of a normal form of the Lorenz system near a triple-zero singularity.
The combination of analytical and numerical tools allows to obtain partial interesting information on the 
complicated dynamics exhibited by system  (\ref{unfoldingFNLorenz}) 
related to a Takens--Bogdanov bifurcation and 
a diagonalizable double-zero degeneracy.
Specifically,  a degenerate heteroclinic connection,  among other global connections, 
gives rise to infinite homoclinic orbits that will lead to the existence of chaos. 
Completing in the future the analysis of system (\ref{unfoldingFNLorenz}) will shed new light on the behavior of the Lorenz system near its triple-zero singularity.

%
%


\begin{thebibliography}{99}
%
\bibitem{Lo:63}
Lorenz E.N. (1963)
Deterministic non-periodic flows. 
{\em J. Atmos. Sci.} {\bf 20} 130--141.




\bibitem{DoKrOs:11}
Doedel E.J., Krauskopf B., Osinga H.M.  (2011)
Global invariant manifolds in the transition to preturbulence in the Lorenz system.
{\em Indagationes Math.} {\bf 22} 222--240.

\bibitem{BaBlSe:11}
Barrio R., Blesa F., Serrano S. (2011)
Global organization of spiral structures in biparameter space of dissipative
  systems with Shilnikov saddle-foci.
{\em Phys. Rev. E} {\bf 84} 035201.

\bibitem{BaShSh:12}
Barrio R., Shilnikov A.L.,  Shilnikov L. (2012)
Kneadings, symbolic dynamics and painting Lorenz chaos.
{\em Int. J. Bifurcation Chaos} {\bf 22} 1230016.

\bibitem{AlFeMeRo:14b}
Algaba A., Fern\'andez-S\'anchez F., Merino M., Rodr\'{\i}guez-Luis A.J. (2014)
Centers on center manifolds in the Lorenz, Chen and L\"u systems. 
{\em Commun. Nonlinear Sci. Numer. Simul.} {\bf 19} 772--775.

\bibitem{AlFeMeRo:15}
Algaba A., Fern\'andez-S\'anchez F., Merino M., Rodr\'{\i}guez-Luis A.J. (2015)
Analysis of the T-point-Hopf bifurcation in the Lorenz system.
{\em Commun. Nonlinear Sci. Numer. Simul.} {\bf 22} 676--691.

\bibitem{CrKrOs:15a}
Creaser J.L., Krauskopf B., Osinga H.M. (2015)
$\alpha$-flips and T-points in the Lorenz system.
{\em Nonlinearity} {\bf 28} R39--R65.

\bibitem{DoKrOs:15}
Doedel E.J., Krauskopf B., Osinga H.M. (2015)
Global organization of phase space in the transition to chaos in the Lorenz system.
{\em Nonlinearity} {\bf 28} R113--R139.

\bibitem{AlDoMeRo:15}
Algaba A., Dom\'{\i}nguez-Moreno M.C., Merino M., Rodr\'{\i}guez-Luis A.J. (2015)
Study of the Hopf bifurcation in the Lorenz, Chen and L\"u systems.
{\em Nonlinear Dynam.} {\bf 79} 885--902.

\bibitem{AlDoMeRo:16}
Algaba A., Dom\'{\i}nguez-Moreno M.C., Merino M., Rodr\'{\i}guez-Luis A.J. (2016)
Takens--Bogdanov bifurcations of equilibria and periodic orbits in the Lorenz system.
{\em Commun. Nonlinear Sci. Numer. Simul.} {\bf 30} 328--343.

\bibitem{AlGaMeRo:16}
Algaba A., Gamero E., Merino M., Rodr\'{\i}guez-Luis A.J. (2016)
Resonances of periodic orbits in the Lorenz system.
{\em Nonlinear Dynam.} {\bf 84} 2111--2136.

\bibitem{AlMeRo:16}
Algaba A., Merino M., Rodr\'{\i}guez-Luis A.J. (2016)
Superluminal periodic orbits in the Lorenz system.
{\em Commun. Nonlinear Sci. Numer. Simul.} {\bf 39} 220--232.

\bibitem{CrKrOs:17}
Creaser J.L., Krauskopf B., Osinga H.M. (2017)
Finding first foliation tangencies in the Lorenz system.
{\em SIAM J. Appl. Dyn. Syst.} {\bf 16} 2127--2164.

\bibitem{Os:18}
Osinga H.M. (2018)
Understanding the geometry of dynamics: the stable manifold of the Lorenz system, 
J. Roy. Soc. New Zeal. {\bf 48} 203--214.

\bibitem{AlDoMeRo:18}
Algaba A., Dom\'{\i}nguez-Moreno M.C., Merino M., Rodr\'{\i}guez-Luis A.J. (2018) 
A Review on Some Bifurcations in the Lorenz System. 
In: Carmona V., Cuevas-Maraver J., Fern\'andez-S\'anchez F., Garc\'ia-Medina E. (eds.) 
Nonlinear Systems, Vol. 1. Understanding Complex Systems. Springer, Cham. 

\bibitem{AlDoMeRo:20}
Algaba A., Dom\'{\i}nguez-Moreno M.C., Merino M., Rodr\'{\i}guez-Luis A.J. (2020)
Double-zero degeneracy and heteroclinic cycles in a perturbation of the Lorenz system.
Preprint.


\bibitem{Shim:80}
Shimizu T., Morioka N. (1980)
On the bifurcation of a symmetric limit cycle to an asymmetric one in a simple model.
{\em Phys. Lett. A}  {\bf 76} 201--204.

\bibitem{Shil:93}
Shil'nikov A.L. (1993)
On bifurcations of the Lorenz attractor in the Shimizu-Morioka model.
{\em Physica D}  {\bf 62} 338--346.

\bibitem{Ru:93}
Rucklidge A.M. (1993)
Chaos  in a low-order model of magnetoconvection.
{\em Physica D}  {\bf 62} 323--337.
 
\bibitem{Liu:04}
Liu C., Liu T., Liu L., Liu K. (2004)
A new chaotic attractor.
{\em Chaos Soliton Fract.} {\bf 22} 1031--1038.

\bibitem{Me:08}
Mello L.F., Messias M., Braga D.C. (2008)
Bifurcation analysis of a new Lorenz-like chaotic system.
{\em Chaos Soliton Fract.} {\bf 37} 1224--1255.
  
\bibitem{KoRo:04}
Kokubu H., Roussarie R. (2004) 
Existence of a singulary degenerate heteroclinic cycle in the Lorenz system
and its dynamical consequences: Part I. 
{\em J. Dyn. Differ. Equ.} {\bf 16} 513--557.
 
\bibitem{GuHo:83}
Guckenheimer, J., Holmes, P.J. (1983)
Nonlinear Oscillations, Dynamical Systems, and Bifurcations of Vector Fields. 
Springer, New York.

\bibitem{Wi:03}
Wiggins S. (2003)
Introduction to Applied Dynamical Systems and Chaos. 
Springer, New York.

\bibitem{Ku:04}
Kuznetsov, Y.A. (2004) 
Elements of Applied Bifurcation Theory. 
Springer, New York.

\bibitem{LiRo:90}
Li C., Rousseau C. (1990)
Codimension 2 symmetric homoclinic bifurcations and application to 1:2 resonance.
{\em Can. J. Math.} {\bf 42} 191--212.

\bibitem{RoFrPo:91}
Rodr\'iguez-Luis A.J., Freire E., Ponce, E. (1991)
On a codimension 3 bifurcation arising in an autonomous electronic circuit, 
in Bifurcation and Chaos: Analysis, Algorithms, Applications, R. Seydel et al. (eds.), 
International Series of Numerical Mathematics, vol. 97, pp. 301--306, Birkh\"auser, Basel.
 
\bibitem{DoOl:12}
Doedel E.J. et al. (2012)
AUTO-07P: Continuation and bifurcation software for
ordinary differential equations. Technical report,
Concordia University.


\bibitem{GlSp:86}
Glendinning P., Sparrow C. (1986)
T-points: a codimension two heteroclinic bifurcation, 
{\em J. Statist. Phys.} {\bf 43} 479--488.


\bibitem{FeFrRo:02} 
Fern\'andez--S\'anchez F., Freire E., Rodr\'{\i}guez-Luis A.J. (2002)
T-Points in a $\mathbb{Z}_2$-symmetric electronic oscillator.
{\em Nonlinear Dynam.} {\bf 28} 53--69.

\bibitem{FeFrRo:08} 
Fern\'andez--S\'anchez F., Freire E., Rodr\'{\i}guez-Luis A.J. (2008)
Analysis of the T-point--Hopf bifurcation.
{\em Physica D} {\bf 237} 292--305.


\end{thebibliography}
\end{document}